# CYCLES ALGÉBRIQUES
# SUR LES SURFACES K3 RÉELLES

Frédéric Mangolte

Dipartimento di Matematica, Via Buonarroti, 2, 56127 PISA, ITALIE
Tel. : (39) 50 59 95 08,   Fax : (39) 50 59 95 24
e-mail : mangolte@dm.unipi.it

Résumé. Sur une surface K3 algébrique réelle $X(\mathbb{R})$, on donne toutes les valeurs possibles pour la dimension $h^1_{\mathrm{alg}}(X(\mathbb{R}))$ du groupe $H^1_{\mathrm{alg}}(X(\mathbb{R}), \mathbb{Z}/2)$ des cycles algébriques de $X(\mathbb{R})$. En particulier, on montre que si $X(\mathbb{R})$ n'est pas une M-surface, on peut toujours déformer $X(\mathbb{R})$ en $X'(\mathbb{R})$ avec $h^1_{\mathrm{alg}}(X'(\mathbb{R})) = \dim H^1(X(\mathbb{R}), \mathbb{Z}/2)$.

On en déduit que dans certains espaces des modules de surfaces K3 algébriques réelles, les classes d'isomorphie réelles de surfaces K3 telles que $h^1_{\mathrm{alg}}(X(\mathbb{R})) \geq k$ forment une réunion dénombrable de sous-espaces de dimension $20 - k$.

## Algebraic cycles on real K3 surfaces

Abstract. For a real algebraic K3 surface $X(\mathbb{R})$, we give all possible values of the dimension $h^1_{\mathrm{alg}}(X(\mathbb{R}))$ of the group $H^1_{\mathrm{alg}}(X(\mathbb{R}), \mathbb{Z}/2)$ of algebraic cycles of $X(\mathbb{R})$. In particular, we prove that if $X(\mathbb{R})$ is not an M-surface, $X(\mathbb{R})$ can always be deformed to some $X'(\mathbb{R})$ with $h^1_{\mathrm{alg}}(X'(\mathbb{R})) = \dim H^1(X(\mathbb{R}), \mathbb{Z}/2)$.

Furthermore, we obtain that in certain moduli space of real algebraic K3 surfaces, the collection of real isomorphism classes of K3 surfaces $X(\mathbb{R})$ such that $h^1_{\mathrm{alg}}(X(\mathbb{R})) \geq k$ is a countable union of subspaces of dimension $20 - k$.

## 1. Introduction

On peut interpréter le nombre de Picard d'une surface algébrique comme le rang du groupe des cycles algébriques modulo équivalence homologique.

Dans l'espace des modules des surfaces K3 algébriques, les classes d'isomorphie de surfaces ayant un nombre de Picard $\rho \geq k$ forment une réunion dénombrable de sous-espaces de dimension $20 - k$ (cf. e.g. [GrHa], p. 594). L'un des objectifs de cet article est de démontrer une propriété similaire pour les surfaces K3 algébriques réelles.

1991 *Mathematics Subject Classification*. 14C25 14P99 14J28.
Supported by HCM contract n° ERBCHBG CT 920011 (MAP Project)

Soit $X(\mathbb{R})$ l'ensemble des points réels d'une variété algébrique $X$ définie sur $\mathbb{R}$, on considère $Y \subset X(\mathbb{R})$ l'ensemble des points réels d'une sous-variété algébrique irréductible définie sur $\mathbb{R}$ de codimension 1 dans $X$. Si $Y$ est de codimension 1 dans $X(\mathbb{R})$, on peut lui associer sa classe fondamentale $[Y] \in \mathrm{H}_{n-1}(X(\mathbb{R}), \mathbb{Z}/2)$ où $n$ est la dimension de $X$. On note $\mathrm{H}^1_{\mathrm{alg}}(X(\mathbb{R}), \mathbb{Z}/2)$ le sous-groupe de $\mathrm{H}^1(X(\mathbb{R}), \mathbb{Z}/2)$ engendré par les classes duales de Poincaré de telles classes fondamentales $[Y]$. C'est le groupe des cycles algébriques de $X(\mathbb{R})$, cf. e.g. [BK]. On note $h^1_{\mathrm{alg}}(X(\mathbb{R})) = \dim \mathrm{H}^1_{\mathrm{alg}}(X(\mathbb{R}), \mathbb{Z}/2)$.

Dans l'espace des modules des surfaces K3 algébriques réelles vérifiant certaines conditions qui seront explicitées plus loin, les classes d'isomorphie réelles de surfaces K3 telles que $h^1_{\mathrm{alg}}(X(\mathbb{R})) \geq k$ forment une réunion dénombrable de sous-espaces de dimension $20 - k$, corollaire 5.8.

On obtient aussi des résultats valables pour les surfaces K3 qui ne sont pas algébriques, on utilisera donc une notion plus large de surface réelle :

Soit $X(\mathbb{C})$ une surface complexe compacte kählerienne, on suppose $X(\mathbb{C})$ munie d'une involution anti-holomorphe $\sigma$, on dit que $(X, \sigma)$ est une surface réelle. On note $X(\mathbb{R})$ l'ensemble des points fixes de $\sigma \colon X(\mathbb{C}) \to X(\mathbb{C})$, c'est la partie réelle de $X$. Soit $Y(\mathbb{C}) \subset X(\mathbb{C})$ une sous-variété analytique irréductible de codimension 1 invariante par $\sigma$, on note $Y(\mathbb{R})$ la partie fixe de $Y(\mathbb{C})$ par $\sigma$. Si $Y(\mathbb{R})$ est de codimension 1 dans $X(\mathbb{R})$, on peut lui associer sa classe fondamentale $[Y(\mathbb{R})] \in \mathrm{H}_1(X(\mathbb{C}), \mathbb{Z}/2)$, cf. [BoHa]. On note $\eta(Y) \in \mathrm{H}^1_{\mathrm{an}}(X(\mathbb{R}), \mathbb{Z}/2)$ la classe duale de Poincaré de $[Y(\mathbb{R})]$, $\eta(Y) = 0$ si $Y(\mathbb{R})$ n'est pas de codimension 1 dans $X(\mathbb{R})$. On note $\mathrm{H}^1_{\mathrm{an}}(X(\mathbb{R}), \mathbb{Z}/2)$ le sous-groupe de $\mathrm{H}^1(X(\mathbb{R}), \mathbb{Z}/2)$ engendré par les classes de cohomologie de la forme $\eta(Y)$, on note $h^1_{\mathrm{an}}(X(\mathbb{R})) = \dim \mathrm{H}^1_{\mathrm{an}}(X(\mathbb{R}), \mathbb{Z}/2)$.

*Remarque.* La notion de cycle analytique de $X(\mathbb{R})$ est différente de celle de cycle analytique réel.

L'application $\eta$ définit donc une correspondance (que l'on précisera plus loin) entre cycles analytiques de $X(\mathbb{C})$ et cycles analytiques de $X(\mathbb{R})$ mais cette correspondance n'est pas triviale :

*(1.1) Exemple.* Soit la surface $S^1 \times S^1$ réalisée comme quartique de $\mathbb{P}^3(\mathbb{R})$, cf. définition 2.6, en changeant légèrement les paramètres on peut en faire une surface lisse et générique. C'est une surface K3 dont la partie réelle $X(\mathbb{R})$ est un tore. La section hyperplane $Y(\mathbb{C})$ nous donne un cycle algébrique non trivial dans $\mathrm{H}^2(X(\mathbb{C}), \mathbb{Z})$ mais aucun cycle algébrique dans la partie réelle : si l'intersection du plan et du tore n'est pas vide, elle est constituée d'une réunion de deux cercles homologues ou d'un cercle homologue à 0 donc $\eta(Y) = 0$.

Pour tenir compte de ce phénomène, on a défini dans [Ma], un morphisme canonique entre la cohomologie de $X(\mathbb{C})$ et la cohomologie de $X(\mathbb{R})$, cf. théorème 2.4.

Comme première conséquence, on trouve des restrictions pour les surfaces de $\mathbb{P}^3(\mathbb{R})$ à être contractibles en un point cf. lemme 2.7 et proposition 2.10.

Après avoir rappelé certains résultats connus sur la topologie des surfaces K3 réelle, on établit un résultat d'obstruction pour une surfaces réelle à vérifier

$$\mathrm{H}^1_{\mathrm{an}}(X(\mathbb{R}), \mathbb{Z}/2) = \mathrm{H}^1(X(\mathbb{R}), \mathbb{Z}/2) .$$

En particulier, une M-surface ne peut vérifier cette égalité, cf. proposition 3.2, corollaire 3.4. On montre ensuite que pour les surfaces K3, c'est la seule obstruction, théorème 4.4. Le corollaire 5.8 annoncé plus haut pour les surfaces K3 algébriques est alors prouvé pour les surfaces K3 générales, corollaire 4.11.

## 2. Topologie des surfaces réelles

Le groupe d'ordre 2 engendré par $\sigma$ est noté $G$. Soit $A$ un $\mathbb{Z}$-module libre muni d'une action de $G$, on note $A^G = \ker(1-\sigma)A$ le sous-groupe invariant et $A(1)^G$ le sous-groupe formé des éléments anti-invariants. On note $\mathrm{H}^i(G, A)$ le i-ème groupe de cohomologie de Galois de $A$. C'est un $\mathbb{Z}/2$-espace vectoriel et on a

$$\mathrm{H}^1(G, A) = A(1)^G/(1-\sigma)A \ ;$$
$$\mathrm{H}^2(G, A) = A^G/(1+\sigma)A \ ;$$
$$\dim \mathrm{H}^1(G, A) = \mathrm{rg}\, A - b - \lambda \ ;$$
$$\dim \mathrm{H}^2(G, A) = b - \lambda \ ;$$

où $b = \mathrm{rg}\, A^G$ et $\lambda = \dim_{\mathbb{Z}/2}(1-\sigma)A \otimes \mathbb{Z}/2$ est la caractéristique de Comessatti de $A$.

Soit $(X, \sigma)$ une surface réelle, on note $(b, \lambda)$ les invariants obtenus en considérant le groupe $A = \mathrm{H}^2(X(\mathbb{C}), \mathbb{Z})$. On note pour $\mathbb{K} = \mathbb{R}$ ou $\mathbb{C}$

$$h^i(X(\mathbb{K})) = \dim \mathrm{H}^i(X(\mathbb{K}), \mathbb{Z}/2) \ ;$$
$$h^*(X(\mathbb{K})) = \sum h^i(X(\mathbb{K})) \ .$$

De la théorie de Smith, on tire que $h^*(X(\mathbb{R})) \leq h^*(X(\mathbb{C}))$ et la différence est paire (cf. [Bre]). On rappelle la définition :

**(2.1) Définition.** *Une M-surface (M pour maximale) est une surface réelle $(X, \sigma)$ telle que $h^*(X(\mathbb{R})) = h^*(X(\mathbb{C}))$. Plus généralement, une $(M{-}a)$-surface est une surface réelle telle que*
$$a = \frac{h^*(X(\mathbb{C})) - h^*(X(\mathbb{R}))}{2} \ .$$

**(2.2) Théorème.** *(Krasnov-Silhol) Soit $(X, \sigma)$ une surface réelle lisse telle que $\mathrm{H}_1(X(\mathbb{C}), \mathbb{Z}/2) = 0$ et $X(\mathbb{R}) \neq \emptyset$, on a*
$$\# X(\mathbb{R}) = \frac{2 + \dim \mathrm{H}^2(G, \mathrm{H}^2(X(\mathbb{C}), \mathbb{Z}))}{2}$$
$$h^1(X(\mathbb{R})) = \dim \mathrm{H}^1(G, \mathrm{H}^2(X(\mathbb{C}), \mathbb{Z}))$$

*Preuve.* En effet, sous les hypothèses du théorème 2.2, on sait que $(X, \sigma)$ est Galois-Maximale, cf. Krasnov, [Kr]. Par ailleurs $\mathrm{H}_1(X(\mathbb{C}), \mathbb{Z}/2) = 0$ implique grâce à la suite des coefficients universel (cf. e.g. [Gre]) que l'homologie de $X(\mathbb{C})$ est sans 2-torsion, et dans ce cas, une surface Galois-Maximale vérifie les égalités voulues, cf. [Si], chap. I.

On a un morphisme bien connu (cf. [BoHa], [BCR], [Si], chap. III), construit à partir de $\eta$ (cf. introduction). Soit $(X, \sigma)$ une surface réelle telle que $X(\mathbb{R}) \neq \emptyset$, pour tout $d \in \mathrm{Pic}(X)^G$, on peut trouver un diviseur invariant $D$ qui représente $d$, cf. [Si], I.(4.5). Soit $D = \sum n_i D_i$ sa décomposition en composantes irréductibles, si $D_i$ est invariant par $\sigma$, on a défini $\eta(D_i) \in \mathrm{H}^1_{\mathrm{an}}(X(\mathbb{R}), \mathbb{Z}/2)$, si $D_i$ n'est pas invariant par $\sigma$, on prolonge $\eta$ en posant $\eta(D_i) = 0$. Finalement on pose $\eta(D) = \sum n_i \eta(D_i)$. On peut montrer que si $D$ est de classe nulle dans $\mathrm{Pic}(X)$, alors $\eta(D) = 0$, cf. [Si], III.(1.7). On définit $\alpha \colon \mathrm{Pic}(X)^G \to \mathrm{H}^1(X(\mathbb{R}), \mathbb{Z}/2)$ par

(2.3) $$\alpha(d) = \eta(D)$$

où $D$ est un représentant invariant de $d$. De plus, on peut montrer que $\alpha$ est surjectif sur $\mathrm{H}^1_{\mathrm{an}}(X(\mathbb{R}), \mathbb{Z}/2)$. Mais il n'est pas injectif en général même si $X(\mathbb{C})$ est simplement connexe, cf. exemple 1.1.

Dans [Ma], chap. I, on a construit un morphisme qui généralise $\alpha$ :

**(2.4) Théorème.** *Soit $(X, \sigma)$ une surface réelle telle que $\mathrm{H}_1(X(\mathbb{C}), \mathbb{Z}) = 0$ et $X(\mathbb{R}) \neq \emptyset$, il existe un morphisme canonique surjectif*
$$\varphi_X \colon \mathrm{H}^2(X(\mathbb{C}), \mathbb{Z})(1)^G \to \mathrm{H}^1(X(\mathbb{R}), \mathbb{Z}/2)$$

*qui vérifie :*

   (1) $\ker \varphi_X = (1 - \sigma) \mathrm{H}^2(X(\mathbb{C}), \mathbb{Z})$,

(2) $\forall \gamma, \gamma' \in \mathrm{H}^2(X(\mathbb{C}), \mathbb{Z})(1)^G$, $(\varphi_X(\gamma), \varphi_X(\gamma')) \equiv Q(\gamma, \gamma') \mod 2$, où $(\ ,\ )$ et $Q$ sont les formes induites par les cup-produits sur $\mathrm{H}^1(X(\mathbb{R}), \mathbb{Z}/2)$ et $\mathrm{H}^2(X(\mathbb{C}), \mathbb{Z})(1)^G$.

Cf. [Ma], I.(5.18).

**(2.5) Théorème.** *Soit* $(X, \sigma)$ *une surface réelle telle que* $\mathrm{H}_1(X(\mathbb{C}), \mathbb{Z}) = 0$ *et* $X(\mathbb{R}) \neq \emptyset$, *alors le diagramme suivant est commutatif*

$$\begin{array}{ccc} \mathrm{Pic}(X)^G & \xrightarrow{c_1} & \mathrm{H}^2(X(\mathbb{C}), \mathbb{Z})(1)^G \\ \alpha \downarrow & & \varphi_X \downarrow \\ \mathrm{H}^1_{\mathrm{an}}(X(\mathbb{R}), \mathbb{Z}/2) & \xrightarrow{j} & \mathrm{H}^1(X(\mathbb{R}), \mathbb{Z}/2) \end{array}$$

*où $j$ est l'inclusion.*

Cf. [Ma], II.(3.1).

**(2.6) Définition.** *On dira que $X(\mathbb{R})$ est une surface de $\mathbb{P}^3(\mathbb{R})$ si $(X, \sigma)$ est une surface réelle qui admet un plongement réel*

$$j_{\mathbb{C}} : X(\mathbb{C}) \hookrightarrow \mathbb{P}^3(\mathbb{C})$$

*réel au sens où $j_{\mathbb{C}} \circ \sigma = conj \circ j_{\mathbb{C}}$ (conj est la conjugaison complexe sur $\mathbb{P}^3(\mathbb{C})$).*

Soit $X(\mathbb{R})$ une surface de $\mathbb{P}^3(\mathbb{R})$ et $h_{\mathbb{C}}$ la classe d'un hyperplan dans $\mathrm{H}^2(\mathbb{P}^3(\mathbb{C}), \mathbb{Z})$, on note $l = j_{\mathbb{C}}^*(h_{\mathbb{C}})$, c'est la classe de la section hyperplane dans $L = \mathrm{H}^2(X(\mathbb{C}), \mathbb{Z})$. Clairement $\sigma(l) = -l$ on a donc $l \in L(1)^G$.

**(2.7) Lemme.** *Une surface de $\mathbb{P}^3(\mathbb{R})$ est contractible en un point si et seulement si $l \in (1 - \sigma)L$.*

*Preuve.* D'après le théorème de Lefschetz sur les sections hyperplanes on a pour les groupes d'homotopie de $X(\mathbb{C})$ (cf. [BPV], p. 46) :

$$(2.8) \quad \begin{aligned} \pi_0(X(\mathbb{C}), \mathbb{Z}) &\cong \pi_0(\mathbb{P}^3(\mathbb{C}), \mathbb{Z}), \\ \pi_1(\mathbb{P}^3(\mathbb{C}), \mathbb{Z}) &\to \pi_1(X(\mathbb{C}), \mathbb{Z}) \to 0. \end{aligned}$$

On en déduit que $X(\mathbb{C})$ est connexe et simplement connexe. En particulier, $\mathrm{H}_1(X(\mathbb{C}), \mathbb{Z}) = 0$ et si $X(\mathbb{R})$ est non vide, on est dans les conditions du théorème 2.4.

La surface $j_{\mathbb{R}}(X(\mathbb{R}))$ est contractible en un point si et seulement si sa classe d'immersion est nulle car $\mathbb{R}^3$ est contractile, cf. [BrWo], p. 85. Cette classe est égale à $j_{\mathbb{R}}^*(h_{\mathbb{R}})$ où $h_{\mathbb{R}}$ est le générateur de $\mathrm{H}^1(\mathbb{P}^3(\mathbb{R}), \mathbb{Z}/2) \cong \mathbb{Z}/2$. D'après le théorème 2.5, on a

$$j_{\mathbb{R}}^*(h_{\mathbb{R}}) = \varphi_X(l).$$

D'après le théorème 2.4, $\varphi_X(l) = 0$ si et seulement si $l \in (1 - \sigma)L$.

**(2.9) Lemme.** *Soit $X(\mathbb{C})$ une surface lisse de $\mathbb{P}^3(\mathbb{C})$, la section hyperplane $l$ n'est pas 2-divisible.*

Cf. [Kh1], lemma 1.2. Pour une preuve plus élémentaire, voir [Ma], I.(6.7).

**(2.10) Proposition.** *Une surface $X(\mathbb{R})$ de $\mathbb{P}^3(\mathbb{R})$ de degré $d = \deg(j_{\mathbb{C}}(X(\mathbb{C})))$ qui vérifie l'une des trois conditions suivantes n'est pas contractible en un point et en particulier n'est pas isotope à une surface de $\mathbb{R}^3 \subset \mathbb{P}^3(\mathbb{R})$.*

  (i) $d \equiv 1 \mod 2$,
  (ii) $(X, \sigma)$ *est une M-surface*,
  (iii) $(X, \sigma)$ *est une (M−1)-surface telle que $d \equiv 0 \mod 4$.*

*Preuve.* Pour simplifier les notations, on note $x \mapsto x^2$ la forme quadratique associée à $Q$, i. e. $x^2 = Q(x,x)$.

  (i) On sait que la partie réelle d'une hypersurface $X(\mathbb{C})$ de $\mathbb{P}^3(\mathbb{C})$ de degré impair est non orientable et ne peut donc être plongée dans $\mathbb{R}^3$.
  (ii) Si $(X, \sigma)$ est une M-surface, $L$ se décompose en somme directe $L = L^G \oplus L(1)^G$ et on a $(1-\sigma)L = 2L(1)^G$. Si $X(\mathbb{R})$ est contractible en un point, on a $l \in (1-\sigma)L$ d'après le lemme 2.7, ce qui est impossible d'après le lemme 2.9.
  (iii) $X$ étant une (M−1)-surface, on peut trouver une base

$$(a_1, \ldots, a_{b-1}, e, e', f_1, \ldots, f_c)$$

de $L$ (où $c = \text{rg}\, L - b - 1$) telle que $\sigma(a_j) = a_j$, $\sigma(e) = e'$ et $\sigma(f_j) = -f_j$, cf. [Si], p. 14. On a alors $(1-\sigma)L = \mathbb{Z}\{(e-e'), 2f_1, \ldots, 2f_c\}$. Si de plus $l \in (1-\sigma)L$, on a

$$(*) \qquad l = n(e-e') + 2(m_1 f_1 + \cdots + m_c f_c) \, ;$$

$$(**) \qquad l^2 \equiv n^2(e-e')^2 \mod 4 \, .$$

Calculons $(e-e')^2$. On note $F = \mathbb{Z}\{f_1, \ldots, f_c\}$. En utilisant la propriété $Q(d, d') = Q(\sigma(d), \sigma(d'))$, on obtient pour tout $i$,

$$Q(e-e', f_i) \equiv 0 \mod 2 \, .$$

De là, le déterminant de $Q_{|L(1)^G}$ vérifie $\det Q_{|L(1)^G} \equiv (e-e')^2 \det Q_{|F}$ mod 4. Par ailleurs, $\det Q_{|L(1)^G} = \pm 2$ ([Si], p. 28). D'où $(e-e')^2 \not\equiv 0$ mod 4.

Maintenant, si $l$ est une section hyperplane, $n$ est impair car sinon, $l$ serait 2-divisible d'après $(*)$ ce qui est impossible d'après le lemme 2.9. Et si $n$ est impair, le degré de $X$ n'est pas divisible par 4 d'après $(**)$.

*(2.11) Remarques.*

(1) On a des exemples de (M−2)-surfaces quartiques de $\mathbb{P}^3(\mathbb{R})$ contractibles en un point, cf. [Kh2] et [Kh3].

(2) Si le degré est pair mais non multiple de 4, il existe des (M−1)-surfaces contractibles en un point. Par exemple, en degré 2, on a la quadrique d'équation $x^2 + y^2 + z^2 - w^2 = 0$ dont la partie réelle est une sphère d'où $h^*(X(\mathbb{R})) = 2$. D'autre part $h^*(X(\mathbb{C})) = 4$ car $X(\mathbb{C}) = \mathbb{P}^1(\mathbb{C}) \times \mathbb{P}^1(\mathbb{C})$ donc c'est une (M−1)-surface.

Une surface K3 est une surface complexe complète $X$ de classe canonique triviale et de premier nombre de Betti nul. D'après un résultat de Siu, on sait que toute surface K3 est kählerienne, cf. [X], exposé XII. On sait que $X(\mathbb{C})$ est simplement connexe, donc $H_1(X(\mathbb{C}), \mathbb{Z}) = 0$ et l'homologie de $X(\mathbb{C})$ est sans torsion. De plus, on a $\operatorname{rg} H^2(X(\mathbb{C}), \mathbb{Z}) = 22$.

Dans toute la suite, $(X, \sigma)$ désigne une surface K3 réelle telle que $X(\mathbb{R}) \neq \emptyset$. D'après le théorème 2.2, on a

$$\begin{aligned} \#X(\mathbb{R}) &= \frac{2 + b - \lambda}{2} \\ h^1(X(\mathbb{R})) &= 22 - b - \lambda \end{aligned} \quad (2.12)$$

et donc

$$\begin{aligned} h^*(X(\mathbb{R}), \mathbb{Z}/2) &= 24 - 2\lambda \\ \chi(X(\mathbb{R})) &= 2b - 20 \,. \end{aligned} \quad (2.13)$$

En considérant les restrictions imposées aux valeurs possibles de $(b, \lambda)$, cf. [Si], en particulier $\lambda \leq \min(b, 22 - b)$, les égalités 2.13 donnent 64 cas pour les valeurs de $(h^*(X(\mathbb{R})), \chi(X(\mathbb{R})))$. Elles sont récapitulées par le diagramme suivant :

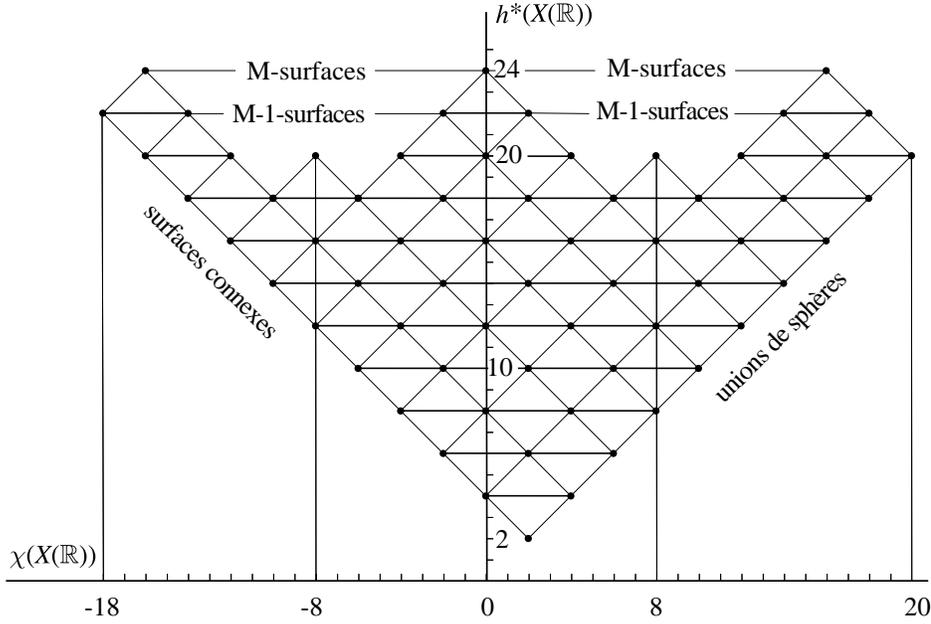

(2.14) Valeurs possibles de $(h^*(X(\mathbb{R})), \chi(X(\mathbb{R})))$ pour une surface K3 réelle.

On peut déduire de ce diagramme et de la proposition suivante qu'il y a 66 types topologiques pour les parties réelles $X(\mathbb{R})$ de surfaces K3 réelles (y compris le type $X(\mathbb{R}) = \emptyset$).

**(2.15) Proposition.** *(Kharlamov) Soit $(X, \sigma)$ une surface K3 réelle, si $X(\mathbb{R})$ est non vide alors $X(\mathbb{R})$ est orientable et comprends au plus une composante connexe non sphérique sauf si $X(\mathbb{R})$ est réunion de deux tores.*

**(2.16) Théorème.** *(Kharlamov) Pour chaque type topologique décrit par le diagramme 2.14 et la proposition 2.15, il existe une surface K3 réelle telle que $X(\mathbb{R}))$ réalise ce type topologique. De plus chaque type topologique peut être réalisé par une quartique de $\mathbb{P}^3(\mathbb{R})$.*

Cf. [Kh2].

## 3. Cycles algébriques et application des périodes

Soit $(X, \sigma)$ une surface réelle telle que $H^1(X(\mathbb{C}), \mathbb{Z}) = 0$, l'application première classe de Chern
$$c_1 : \text{Pic}(X) \hookrightarrow H^2(X(\mathbb{C}), \mathbb{Z})$$
est alors injective. On identifie $H^2(X(\mathbb{C}), \mathbb{Z})$ avec le sous-groupe $H^2(X(\mathbb{C}), \mathbb{Z}) \otimes 1$ de l'espace $H^2(X(\mathbb{C}), \mathbb{C})$.

**(3.1) Théorème.** *(Lefschetz) Soit $(X, \sigma)$ une surface réelle qui vérifie*

$$\mathrm{H}^1(X(\mathbb{C}), \mathbb{Z}) = 0 \ ,$$

*alors $c_1$ induit un isomorphisme*

$$\mathrm{Pic}(X)^G \cong \mathrm{H}^2(X(\mathbb{C}), \mathbb{Z})(1)^G \cap \mathrm{H}^{1,1}(X) \ .$$

*Preuve.* Du fait que la classe de Chern d'un diviseur invariant est anti-invariante, cf. [Si], I.(4.12), ce théorème découle du théorème de Lefschetz sur les (1,1)-cycles, cf. [BPV], p. 119.

**(3.2) Proposition.** *Soit $(X, \sigma)$ une $(M-a)$-surface de genre géométrique $p_g = h^{2,0}$ telle que $\mathrm{H}_1(X(\mathbb{C}), \mathbb{Z}) = 0$, on a*

$$h^1_{\mathrm{an}}(X(\mathbb{R})) \leq h^1(X(\mathbb{R})) - (p_g - a) \ .$$

*En particulier, si $a < p_g$, on a*

$$h^1_{\mathrm{an}}(X(\mathbb{R})) < h^1(X(\mathbb{R})) \ .$$

En remarquant que pour une $(M-a)$-surface qui vérifie $\mathrm{H}_1(X(\mathbb{C}), \mathbb{Z}) = 0$, on a

$$h^1(X(\mathbb{R})) = \mathrm{rg}\, \mathrm{H}^2(X(\mathbb{C}), \mathbb{Z})(1)^G - a \ ,$$

cf. [Si], II.(3.3), la proposition découle du lemme suivant,

**(3.3) Lemme.** *Soit $X$ une surface algébrique réelle telle que $\mathrm{H}^1(X(\mathbb{C}), \mathbb{Z}) = 0$, on a*

$$h^1_{\mathrm{an}}(X(\mathbb{R})) \leq \mathrm{rg}\, \mathrm{H}^2(X(\mathbb{C}), \mathbb{Z})(1)^G - p_g \ .$$

*Preuve du lemme.* Sous ces hypothèse, $c_1$ est injective. Soit $\mathcal{B}$ une base de $\mathrm{H}^1_{\mathrm{an}}(X(\mathbb{R}), \mathbb{Z}/2)$, on peut la remonter par $\alpha$ (cf. 2.3) en une famille libre de $\mathrm{Pic}(X)^G$ dont l'image par $c_1$ est encore une famille libre de cardinal $h^1_{\mathrm{an}}(X(\mathbb{R}))$ que l'on note $\mathcal{B}'$. D'après le théorème 3.1, on a

$$\mathcal{B}' \subset \mathrm{H}^{1,1}(X) \cap \mathrm{H}^2(X(\mathbb{C}), \mathbb{Z})(1)^G \ .$$

Soit $p_g$ formes holomorphes linéairement indépendantes dans $\mathrm{H}^{2,0}(X)$, a chacune de ces formes $\omega$, on associe $y = \mathrm{Im}(\omega + \bar{\omega}^\sigma)$. Alors $y$ est anti-invariant par $\sigma$ et par construction on a

$$y \in (\mathrm{H}^{2,0}(X) \oplus \mathrm{H}^{0,2}(X)) \cap \mathrm{H}^2(X(\mathbb{C}), \mathbb{R})(1)^G \ .$$

Donc $y$ est orthogonal à $\mathcal{B}'$ dans $\mathrm{H}^2(X(\mathbb{C}), \mathbb{R})$. Si l'on note $B'$ le sous-espace de $\mathrm{H}^2(X(\mathbb{C}), \mathbb{R})$ engendré par $\mathcal{B}'$, on a

$$\dim B' \leq \mathrm{rg}\, \mathrm{H}^2(X(\mathbb{C}), \mathbb{Z})(1)^G - p_g \ .$$

**(3.4) Corollaire.** *Soit $X$ une M-surface telle que $p_g > 0$ et $H_1(X(\mathbb{C}), \mathbb{Z}) = 0$, on a alors*

$$h^1_{\mathrm{an}}(X(\mathbb{R})) < h^1(X(\mathbb{R})) \ .$$

*(3.5) Remarque.* Ce résultat reste vrai si $H_1(X(\mathbb{C}), \mathbb{Z}) \neq 0$, mais la preuve est plus longue, cf. [Ma], chap II.

Dans toute la suite, si $M$ désigne une $\mathbb{Z}$-module et $\mathbb{K} = \mathbb{R}$ ou $\mathbb{C}$, on note $M_\mathbb{K}$ le produit tensoriel $M \otimes \mathbb{K}$. De même si $f : M \to M'$ est un morphisme, on note

$$f_\mathbb{K} : M_\mathbb{K} \to M'_\mathbb{K}$$

le morphisme tensorisé par $id_\mathbb{K}$ et si $Q$ est une forme bilinéaire sur $M$, $Q_\mathbb{K}$ désigne l'extension $\mathbb{K}$-bilinéaire de $Q$ à $M_\mathbb{K}$.

On note $L$ un $\mathbb{Z}$-module muni d'une forme bilinéaire entière symétrique non dégénérée paire $Q$ de signature $(3, 19)$. Le couple $(L, Q)$ est unique à isomorphisme près, cf. [Se].

Soit $X$ une surface K3, on sait ([X], exposé IV), que $H^2(X(\mathbb{C}), \mathbb{Z})$ muni de la forme induite par le cup-produit est isométrique à $(L, Q)$. Soit $f : H^2(X(\mathbb{C}), \mathbb{Z}) \to L$ une telle isométrie. On dit que le couple $(X, f)$ est une surface K3 marquée. Soit la décomposition de Hodge

$$H^2(X(\mathbb{C}), \mathbb{C}) = H^{2,0}(X) \bigoplus H^{1,1}(X) \bigoplus H^{0,2}(X) \ ;$$

on a $h^{2,0} = \dim H^{2,0}(X) = 1$ et $h^{1,1} = \dim H^{1,1}(X) = 20$ (cf. [X], exposé III).

Soit $(X, f)$ une surface K3 marquée. On considère $P \subset L_\mathbb{R}$ l'image par $f_\mathbb{R}$ du sous-espace $H^2(X(\mathbb{C}), \mathbb{R}) \cap (H^{2,0}(X) \oplus H^{0,2}(X))$ de $H^2(X(\mathbb{C}), \mathbb{R})$. Du fait que $H^{2,0}(X)$ est de dimension complexe 1, $P$ est un 2-plan de $L_\mathbb{R}$. On choisit une orientation de $P$ de telle sorte que pour toute 2-forme holomorphe $\omega \in H^{2,0}(X)$, la base $(\mathrm{Re}(\omega), \mathrm{Im}(\omega))$ soit directe. Nous dirons que ce plan orienté $P$ est la période de $(X, f)$.

**(3.6) Théorème.** *(Kharlamov) Soit $\sigma$ une involution de $(L, Q)$, on note $b = \mathrm{rg}\, L^G$. Il existe une surface K3 marquée $(X, f)$ telle que $\tilde{\sigma} = f^{-1} \circ \sigma \circ f$ provienne d'une structure réelle sur $X$ si et seulement si la restriction de la forme bilinéaire à $L^G$ a pour signature $(1, b - 1)$.*

Cf. [Si], p. 184.

Soit $\sigma$ une involution sur $L$ satisfaisant la condition du théorème 3.6. Soit $C_1$ le cône dans $L^G_\mathbb{R}$ défini par $x^2 > 0$. Du fait que par hypothèse, la restriction de $Q_\mathbb{R}$ à $L_{\mathbb{R}}^G$ a pour signature $(1, b-1)$, $C_1$ se décompose en 2 composantes connexes $C_1^+$ et $C_1^-$.

On note $\Omega_1^+(\sigma)$ l'image de $C_1^+$ dans $\mathbb{P}(L_{\mathbb{R}}^G) \cong \mathbb{P}^{b-1}(\mathbb{R})$.

Soit $L_{\mathbb{R}}(1)^G = \{l \in L_{\mathbb{R}}/\sigma(l) = -l\}$ et $C_2$ le cône dans $L_{\mathbb{R}}(1)^G$ défini par $y^2 > 0$. L'hypothèse faite sur $\sigma$ implique que la restriction de $Q_{\mathbb{R}}$ à $L_{\mathbb{R}}(1)^G$ a pour signature $(2, 20-b)$, d'où $C_2$ est connexe.

On note $\Omega_2(\sigma)$ l'image de $C_2$ dans l'ensemble des droites orientées de $L_{\mathbb{R}}(1)^G$. On peut considérer $\Omega_2(\sigma)$ comme un ouvert de l'hypersphère $S^{22-b-1}$.

Soit $e \in L^G$ tel que $e^2 = -2$, les hyperplans $H_e = \{x/Q(x,e) = 0\}$ de $\mathbb{P}(L_{\mathbb{R}}^G)$ découpent des ouverts dans $\Omega_1^+(\sigma)$.

On note $\tilde{\Omega}_1^+(\sigma)$ la réunion de ces ouverts.

**(3.7) Théorème.** *(Nikulin-Silhol) Soit $\sigma$ une involution de $(L, Q)$ satisfaisant la condition du théorème 3.6, l'ensemble des classes d'isomorphie de surfaces K3 marquées $(X, f)$ pour lesquelles $\sigma$ provient d'une structure réelle de $X$ est paramétré bijectivement par l'ouvert*

$$\Omega(\sigma) = \tilde{\Omega}_1^+(\sigma) \times \Omega_2(\sigma) \ de \ \mathbb{P}^{b-1}(\mathbb{R}) \times S^{22-b-1}$$

Cf. [Si], p. 186.

**(3.8) Définition.** *Soit $\sigma$ une involution sur $L$ qui vérifie l'hypothèse du théorème 3.6, on appelle surface K3 de type $\sigma$ toute surface K3 réelle marquée $(X, f)$ de période appartenant à $\Omega(\sigma)$, ce que l'on notera abusivement $(X, f) \in \Omega(\sigma)$.*

## 4. Cycles analytiques sur les surfaces K3 kähleriennes réelles

On note
$$r \colon L(1)^G \longrightarrow L(1)^G/(1-\sigma)L = \mathrm{H}^1(G, L)$$

la surjection canonique. Soit $(X, f)$ une surface K3 réelle marquée. On note $\hat{f}$ la restriction de $f$ à $\mathrm{H}^2(X(\mathbb{C}), \mathbb{Z})(1)^G$.

On considère sur $\mathrm{H}^1(G, L)$ la forme à coefficients dans $\mathbb{Z}/2$ induite par $Q$, on considère sur $\mathrm{H}^1(X(\mathbb{R}), \mathbb{Z}/2)$ et $\mathrm{H}^2(X(\mathbb{C}), \mathbb{Z})(1)^G$, les formes induites par le cup-produit.

**(4.1) Définition.** *Soit $\gamma \in \mathrm{H}^1(X(\mathbb{R}), \mathbb{Z}/2)$ et $\tilde{\gamma} \in \varphi_X^{-1}\{\gamma\}$. On pose $\bar{f}(\gamma) = r \circ \hat{f}(\tilde{\gamma})$.*

D'après 2.4, $\bar{f}$ est bien définie et c'est une isométrie qui rend commutatif le diagramme suivant :

$$
(4.2) \quad \begin{array}{ccc} \mathrm{H}^2(X(\mathbb{C}), \mathbb{Z})(1)^G & \xrightarrow[\sim]{\hat{f}} & L(1)^G \\ \varphi_X \downarrow & & \downarrow r \\ \mathrm{H}^1(X(\mathbb{R}), \mathbb{Z}/2) & \xrightarrow[\sim]{\bar{f}} & \mathrm{H}^1(G, L) \end{array}
$$

**(4.3) Définition.** *Soit $(X, \sigma_X)$ une surface réelle, on dit que $(Y, \sigma_Y)$ est une déformation réelle de $(X, \sigma_X)$ s'il existe deux variétés réelles et une submersion holomorphe propre*

$$\pi : (V, \sigma) \to (B, \tau)$$

*qui commute avec $\sigma$ et $\tau$ et qui vérifie la condition suivante :*

*il existe deux points $p_1$ et $p_2$ dans la même composante connexe de $B(\mathbb{R})$ tels que $X(\mathbb{C}) = \pi^{-1}(p_1)$, $Y(\mathbb{C}) = \pi^{-1}(p_2)$ et $\sigma_Z$ ($Z = X$ ou $Y$) est la restriction de $\sigma$ à $Z$.*

**(4.4) Théorème.** *Soit $X$ une surface K3 réelle, si $X$ n'est pas une M-surface, alors pour tout sous-groupe $K$ de $\mathrm{H}^1(X(\mathbb{R}), \mathbb{Z}/2)$, il existe une déformation réelle $Y$ de $X$ et une isométrie*

$$u : \mathrm{H}^1(X(\mathbb{R}), \mathbb{Z}/2) \to \mathrm{H}^1(Y(\mathbb{R}), \mathbb{Z}/2) \text{ telle que :}$$
$$u(K) = \mathrm{H}^1_{\mathrm{an}}(Y(\mathbb{R}), \mathbb{Z}/2) \, .$$

**(4.5) Corollaire.** *Pour tout entier $k$ vérifiant $k \leq h^1(X(\mathbb{R}))$, il existe une déformation réelle $Y$ de $X$ telle que*

$$h^1_{\mathrm{an}}(Y(\mathbb{R})) = k \, .$$

**(4.6) Théorème.** *Soit $X$ une M-surface K3, pour tout entier $k$ vérifiant :*

$$k < h^1(X(\mathbb{R}))$$

*il existe une déformation réelle $Y$ de $X$ telle que*

$$h^1_{\mathrm{an}}(Y(\mathbb{R})) = k \, .$$

*Preuve du théorème 4.4.* Soit $X$ une surface K3 réelle, on peut supposer que $h^1(X(\mathbb{R})) \neq 0$ puisque si $h^1(X(\mathbb{R})) = 0$, l'assertion est trivialement vérifiée.

Soit $f\colon \mathrm{H}^2(X(\mathbb{C}), \mathbb{Z}) \to L$ un marquage de $X$, soit $\omega \in \mathrm{H}^{2,0}(X)$ une forme holomorphe. Alors $\bar{\omega}^\sigma$ est aussi une forme holomorphe, on pose $x = f_\mathbb{R}(\mathrm{Re}(\omega + \bar{\omega}^\sigma))$. En remplaçant $x$ par $-x$ si nécessaire, on peut supposer que $x \in C_1^+$. On pose $t = h^1(X(\mathbb{R}))$, comme $X(\mathbb{R})$ est orientable, on a $t \geq 2$.

Pour un sous-module $N$ de $L(1)^G$, on note $N^\perp$ son orthogonal pris dans $L(1)^G$. On dit que $M \subset N$ est un sous-module primitif si $M = M_\mathbb{R} \cap N$, ou autrement dit si $M$ est facteur direct dans $N$.

**(4.7) Lemme.** *Il existe un sous-module primitif $M \subset L(1)^G$ de rang $t$ tel que $r(M) = \mathrm{H}^1(G, L)$ qui vérifie*
$$C_2 \cap M_\mathbb{R}^\perp \neq \emptyset \ .$$

Ce lemme sera prouvé plus loin. On choisit $y \in C_2 \cap M_\mathbb{R}^\perp$. D'après le théorème 3.7, le plan $P$ défini par $x$ et $y$ et orienté de telle sorte que la base $(x, y)$ soit directe est la période d'une surface K3 marquée $(Y, g)$ sur laquelle $g^{-1} \circ \sigma \circ g$ induit une structure réelle. De plus, comme $C_2$ est connexe, $Y$ est une déformation de $X$, cf. [X], exposé V.

On obtient alors un diagramme similaire au diagramme 4.2 :

$$\begin{array}{ccc} \mathrm{H}^2(Y(\mathbb{C}), \mathbb{Z})(1)^G & \xrightarrow[\sim]{\hat{g}} & L(1)^G \\ {\scriptstyle \varphi_Y}\downarrow & & \downarrow {\scriptstyle r} \\ \mathrm{H}^1(Y(\mathbb{R}), \mathbb{Z}/2) & \xrightarrow[\sim]{\bar{g}} & \mathrm{H}^1(G, L) \end{array}$$

D'autre part, par définition de la période,
$$g_\mathbb{R}^{-1}(P) = (\mathrm{H}^{2,0}(Y) \oplus \mathrm{H}^{0,2}(Y)) \cap \mathrm{H}^2(Y(\mathbb{C}), \mathbb{R}) \ ;$$

on a donc $\hat{g}^{-1}(M) \subset \mathrm{H}^2(Y(\mathbb{C}), \mathbb{Z})(1)^G \cap \mathrm{H}^{1,1}(Y)$ par construction et $\varphi_Y \circ \hat{g}^{-1}(M) = \mathrm{H}^1_{\mathrm{an}}(X(\mathbb{R}), \mathbb{Z}/2)$. On pose $u = \bar{g}^{-1} \circ \bar{f}$ et on a d'après 3.1 et 2.5,
$$u(K) \subset \mathrm{H}^1_{\mathrm{an}}(Y(\mathbb{R}), \mathbb{Z}/2) \ .$$

Si $\dim_{\mathbb{Z}/2} K = t$, on a clairement $u(K) = \mathrm{H}^1_{\mathrm{an}}(Y(\mathbb{R}), \mathbb{Z}/2)$. Si $\dim_{\mathbb{Z}/2} K < t$, les surfaces $Y$ pour lesquelles l'inclusion précédente est une égalité forment un ouvert dense parmi celles qui vérifient l'inclusion. Il est donc loisible de choisir $y \in C_2 \cap M_\mathbb{R}^\perp$ pour que
$$u(K) = \mathrm{H}^1_{\mathrm{an}}(Y(\mathbb{R}), \mathbb{Z}/2) \ .$$

On dit qu'un sous-module $<a, b>$ est un plan hyperbolique si $a$ et $b$ sont isotropes et de produit $Q(a, b) = 1$.

**(4.8) Lemme.** *(Nikulin-Silhol) Soit $\sigma$ une involution de $(L,Q)$ telle que la restriction $Q_{|L^G}$ ait pour signature $(1, b-1)$, le module $L(1)^G$ admet alors une décomposition orthogonale de l'une des trois formes suivantes :*

(1) $L(1)^G = <a_1, b_1> \perp <a_2, b_2> \perp S$,
(2) $L(1)^G = <2> \perp <a_1, b_1> \perp S$,
(3) $L(1)^G = <2> \perp <2> \perp S$

*où $<a_i, b_i>$ est un plan hyperbolique, $Q_{|S}$ est définie négative quand $S \neq 0$ et $<2>$ est un sous-module de $L(1)^G$ engendré par un élément de carré 2.*

*De plus, chaque sous-module du type $<2>$ apparaissant dans ces décompositions est contenu dans $(1-\sigma)L$ et chaque sous-module du type $<a_i, b_i>$ vérifie*

$$<a_i, b_i> \cap (1-\sigma)L \subset 2L \ .$$

Cf. [Si], pp. 190–191 et corollaire VIII.(4.8).

**(4.9) Lemme.** *Soit $N$ un $\mathbb{Z}$-module libre muni d'une forme bilinéaire symétrique non dégénérée $Q$ et $M \subset N$ un sous-module primitif, alors il existe un sous-module primitif $M' \subset N$ tel que*

(1) $M' \equiv M \mod 2N$
(2) *la restriction de $Q$ à $M'$ est non dégénérée.*

*Preuve.* Soit $F = M \cap M^\perp$, toute décomposition $M = F \oplus H$ est orthogonale et $Q_{|H}$ est non dégénérée. Supposons que $F \neq \{0\}$, soit $\gamma \in F$ un élément primitif, $F_1 = \mathbb{Z}\gamma$ et soit $F_2$ un supplémentaire de $F_1$, $F = F_1 \oplus F_2$. Par définition de $F$, $\gamma$ est isotrope et il existe $\gamma' \in H^\perp$ tel que $(\gamma + 2\gamma')^2 \neq 0$. En effet, si pour tout $\gamma' \in H^\perp$, on a $(\gamma + 2\gamma')^2 = 0$, alors pour tout $\gamma' \in H^\perp$, on a $\gamma'^2 = 0$ ce qui est contradictoire puisque $Q$ est non dégénérée. On note $F_1'$ le sous-module de $N$ engendré par $\gamma + 2\gamma'$ et on pose $M' = H \perp F_1' \oplus F_2$. On a alors $M' \equiv M \mod 2N$ et $M' \cap M'^\perp \subset F_2$ donc $\operatorname{rg} M' \cap M'^\perp < \operatorname{rg} F$. Le lemme se prouve par récurence sur le rang de $M' \cap M'^\perp$.

**(4.10) Lemme.** *Soit $N$ un $\mathbb{Z}$-module libre muni d'une forme bilinéaire symétrique non dégénérée $Q$ de signature $(r,s)$ et $M \subset N$ un sous-module primitif tel que la restriction $Q_{|M}$ soit non dégénérée de signature $(r', s')$, on suppose $s - s' > 0$. Si $M$ contient un plan hyperbolique primitif plongé, c'est-à-dire si il existe une décomposition $M = <a, b> \perp S$, il existe un sous-module primitif $M' \subset N$ qui vérifie*

(1) $M' \equiv M \mod 2N$,
(2) *la restriction $Q_{|M'}$ est de signature $(r'-1, s'+1)$ .*

*Preuve.* Comme $s - s' > 0$, il existe $\gamma \in M^\perp$ tel que $\gamma^2 < 0$.

On pose
$$\alpha = a + 2\gamma, \qquad \beta = b + 2\gamma$$

et on a $<a, b> \equiv <\alpha, \beta> \mod 2N$, de plus la restriction $Q_{|<\alpha,\beta>}$ est définie négative.

Le sous-module $M' = <\alpha, \beta> \perp S$ vérifie alors les conditions (1) et (2).

*Preuve du lemme 4.7.*

Du lemme 4.8, on déduit qu'il existe un sous-module primitif $M \subset L(1)^G$ de rang $t$ tel que $r(M) = \mathrm{H}^1(G, L)$ qui admet une décomposition orthogonale de l'une des trois formes suivante

(1) $M = <a_1, b_1> \perp <a_2, b_2> \perp S'$,
(2) $M = <a_1, b_1> \perp S'$,
(3) $M = S'$

de plus, d'après le lemme 4.9, on peut supposer que si $S' \neq 0$, $Q_{|S'}$ est non dégénérée donc définie négative.

Supposons que $C_2 \cap M_{\mathbb{R}}^\perp = \emptyset$. Alors d'après 4.8, $M$ se décompose sous la forme (1). En appliquant le lemme 4.10, on obtient un module $M'$ qui vérifie les conditions du lemme 4.7.

*Preuve du théorème 4.6.* Si $X$ est une M-surface, on a $(1-\sigma)L = 2L$ et $L(1)^G$ se décompose sous la forme (1), cf. lemme 4.8. Il est alors clair qu'il existe un sous-module primitif $M \subset L(1)^G$ de rang $t - 1$ tel que $\dim_{\mathbb{Z}/2} r(M) = t - 1$ et $M$ se décompose sous la forme (1).

Notons $\Omega^k(\sigma)$ le sous-espace de $\Omega(\sigma)$ formé des périodes de surfaces K3 $(X, f)$ de type $\sigma$ telles que $h_{\mathrm{an}}^1(X(\mathbb{R})) \geq k$.

**(4.11) Corollaire.** *Soit $(X, f) \in \Omega(\sigma)$, pour tout entier $k$ vérifiant les conditions de 4.5 ou 4.6, les surfaces $(Y, g)$ de $\Omega(\sigma)$ telles que $h_{\mathrm{an}}^1(Y(\mathbb{R})) \geq k$ forment une réunion dénombrable de sous-variétés de dimension $20 - k$ dans $\Omega(\sigma)$. En particulier,*
$$\dim \Omega^k(\sigma) = 20 - k \ .$$

## 5. Cycles algébriques sur les surfaces K3 algébriques réelles

On considère maintenant le cas où $(X, \sigma)$ est une surface K3 algébrique réelle, c'est-à-dire qu'il existe un diviseur ample $D$ sur $X$ invariant par $\sigma$. Soit $l$ la première classe de Chern de $D$, c'est une polarisation de $X$ et en particulier, $l^2 > 0$. Comme $(X, \sigma)$ est réelle, $D$ est invariant et donc d'après [Si], I.(4.12), la classe $l$ est anti-invariante. On peut définir un espace des périodes adapté à cette situation, voir aussi [X], p. 13 pour le cas complexe,

**(5.1) Définition.** *Soit $\sigma$ une involution de $(L, Q)$ telle que la restriction de $Q$ à $L^G$ ait pour signature $(1, b-1)$ et $l$ un élément primitif de $L(1)^G$ tel que $l^2 > 0$, on note $\Omega_l(\sigma)$ le sous-espace de $\Omega(\sigma)$ formé des couples $(x, y)$ tels que si $\tilde{y}$ appartient à la droite $y$, on ait*

$$Q(\tilde{y}, l) = 0 \ .$$

Le sous-espace de $\Omega_l(\sigma)$ formé des périodes de surfaces K3 $(X, f)$ de type $\sigma$ polarisées par $f^{-1}(l)$ est défini par la condition supplémentaire suivante :

pour tout élément $e \in l^\perp$ tel que $e^2 = -2$, on a $Q(\tilde{y}, e) \neq 0$, cf. [X], exposé X, p. 127.

*(5.2) Remarque.* On a $\dim \Omega_l(\sigma) = 19$ d'après 4.11.

Dans ce cas précis, on a une obstruction supplémentaire. Soit $\sigma$ et $l$ vérifiant les conditions de la définition 5.1, on pose $k_0 = 0$ si $r(l) = 0$ et $k_0 = 1$ si $r(l) \neq 0$.

**(5.3) Proposition.** *Soit $X$ une (M−1)-surface K3 telle que $k_0 = 0$, alors*

$$h^1_{\mathrm{alg}}(X(\mathbb{R})) < h^1(X(\mathbb{R})) \ .$$

*Preuve.* Supposons que $(X, f)$ soit une (M-1)-surface marquée telle que $h^1_{\mathrm{alg}}(X(\mathbb{R})) = h^1(X(\mathbb{R}))$. Si $(x, y)$ est la période de $(X, f)$, on a $l \in y^\perp \cap L(1)^G$. Mais comme $X$ est une (M-1)-surface, on a $\mathrm{rg}(y^\perp \cap L(1)^G) = h^1(X(\mathbb{R}))$ et $r(y^\perp \cap L(1)^G) = \mathrm{H}^1(X(\mathbb{R}), \mathbb{Z}/2)$. Par ailleurs, $l$ est primitif et $M \cap (1 - \sigma)L = 2M$ donc $r(l) \neq 0$.

**(5.4) Théorème.** *Soit $\sigma$ et $l$ vérifiant les conditions de la définition 5.1 et soit $(X, f) \in \Omega_l(\sigma)$, on note $\tilde{l} = f^{-1}(l)$.*

*Si $X$ n'est ni une M-surface ni une (M−1)-surface telle que $k_0 = 0$, alors pour tout sous-groupe $K$ de $\mathrm{H}^1(X(\mathbb{R}), \mathbb{Z}/2)$ contenant $\varphi_X(\tilde{l})$, il existe une déformation réelle $Y$ de $X$ dont l'espace total est inclus dans $\Omega_l(\sigma)$ et une isométrie*

$$u \colon \mathrm{H}^1(X(\mathbb{R}), \mathbb{Z}/2) \to \mathrm{H}^1(Y(\mathbb{R}), \mathbb{Z}/2) \text{ telle que :}$$
$$u(K) = \mathrm{H}^1_{\mathrm{alg}}(Y(\mathbb{R}), \mathbb{Z}/2) \ .$$

**(5.5) Corollaire.** *Pour tout entier $k$ vérifiant $k_0 \leq k \leq h^1(X(\mathbb{R}))$, il existe une déformation réelle $Y$ de $X$ telle que*

$$h^1_{\mathrm{alg}}(Y(\mathbb{R})) = k \ .$$

**(5.6) Théorème.** *Soit $(X, f) \in \Omega_l(\sigma)$, pour tout entier $k$ vérifiant :*

$$\begin{cases} 1 \leq k < h^1(X(\mathbb{R})) & \text{si } X \text{ est une } M\text{-surface}, \\ 0 \leq k < h^1(X(\mathbb{R})) & \text{si } X \text{ est une } (M{-}1)\text{-surface telle que } k_0 = 0, \end{cases}$$

*il existe une déformation réelle $Y$ de $X$ telle que*

$$h^1_{\mathrm{alg}}(Y(\mathbb{R})) = k \ .$$

*Preuve du théorème 5.4.* La preuve est la même que pour le théorème 4.4 en utilisant le lemme suivant en lieu et place du lemme 4.7.

**(5.7) Lemme.** *Il existe un sous-module primitif $M \subset L(1)^G$ de rang $t$ tel que $r(M) = \mathrm{H}^1(G, L)$ qui vérifie*

$$C_2 \cap M_\mathbb{R}^\perp \cap l^\perp \neq \emptyset \ .$$

*Preuve.* Comme dans la preuve du lemme 4.7, on considère un sous-module primitif $M \subset L(1)^G$ de rang $t$ tel que $r(M) = \mathrm{H}^1(G, L)$ qui admet une décomposition orthogonale de l'une des trois formes suivante

(1) $M = <a_1, b_1> \perp <a_2, b_2> \perp S'$,
(2) $M = <a_1, b_1> \perp S'$,
(3) $M = S'$.

Si $k_0 \neq 0$, on considère un élément primitif $l' \in M$ tel que $r(l') = r(l)$ dans $\mathrm{H}^1(G, L)$. Soit un sous-module $R \subset L(1)^G$ tel que $R \oplus \mathbb{Z}l' \equiv M \mod 2L(1)^G$ et qui admet une décomposition orthogonale de la même forme que $M$ ou de l'une des formes :

(4) $R = <a_1> \perp T$,
(5) $R = <a_1> \perp <a_2, b_2> \perp T$.

D'après le lemme 4.9, on peut supposer que si $T \neq 0$, $Q_{|T}$ est non dégénérée donc définie négative. En remplaçant au besoin $a_1$ par $a_1 - 2b_1$ et en utilisant le lemme 4.10, on peut trouver un sous-module $R' \equiv R \mod 2L(1)^G$ tel que $Q_{|R'}$ soit définie négative. Le module $M' = \mathbb{Z}l \oplus R'$ vérifie alors les hypothèses du lemme 5.7.

Si $k_0 = 0$, on raisonne directement avec $M$.

*Preuve du théorème 5.6.* Soit $(X, f) \in \Omega_l(\sigma)$, si $k_0 = 1$, $r(l) \neq 0$ dans $\mathrm{H}^1(G, L)$ et on a

$$\forall\, Y \in \Omega_l(\sigma), \qquad h^1_{\mathrm{alg}}(Y(\mathbb{R})) \geq 1\ .$$

Maintenant si $X$ est une M-surface, on a $k_0 = 1$. En effet, dans ce cas,

$$(1 - \sigma)\,\mathrm{H}^2(X(\mathbb{C}), \mathbb{Z}) = 2\,\mathrm{H}^2(X(\mathbb{C}), \mathbb{Z})$$

et $l$ est primitif.

Notons $\Omega_l^k(\sigma)$ le sous-espace de $\Omega_l(\sigma)$ formé des périodes de surfaces K3 $(X, f)$ telles que $h^1_{\mathrm{alg}}(X(\mathbb{R})) \geq k$.

**(5.8) Corollaire.** *Soit $(X, f) \in \Omega_l(\sigma)$, pour tout entier $k$ vérifiant les conditions de 5.5 ou 5.6, les surfaces $(Y, g)$ de $\Omega_l(\sigma)$ telles que $h^1_{\mathrm{alg}}(Y(\mathbb{R})) \geq k$ forment une réunion dénombrable de sous-variétés de dimension $19 + k_0 - k$ dans $\Omega_l(\sigma)$. En particulier,*

$$\dim \Omega_l^k(\sigma) = 19 + k_0 - k\ .$$

On appelle degré d'une surface K3, $X \in \Omega_l(\sigma)$, l'entier $l^2$. Si $\tilde{l} = f^{-1}(l)$ engendre un plongement $j_{\tilde{l}}$ dans $\mathbb{P}^{\frac{l^2+2}{2}}$, $l^2$ est le degré de $j_{\tilde{l}}(X(\mathbb{R}))$ dans $\mathbb{P}^{\frac{l^2+2}{2}}(\mathbb{R})$.

**(5.9) Lemme.** *Soit $X(\mathbb{R})$ une quartique de $\mathbb{P}^3(\mathbb{R})$, alors $k_0 = 0$ si et seulement si $X(\mathbb{R})$ est contractible en un point.*

*Preuve.* En effet, d'après le lemme 2.9, $\tilde{l}$ est primitif dans $\mathrm{H}^2(X(\mathbb{C}), \mathbb{Z})$ et d'après le lemme 2.7, $\varphi_X(\tilde{l}) = 0$ si et seulement si $X(\mathbb{R})$ est contractible en un point.

**(5.10) Corollaire.** *Si $X(\mathbb{R})$ est une quartique de $\mathbb{P}^3(\mathbb{R})$ et $(X, \sigma)$ une $(M-1)$-surface, alors $k_0 = 1$.*

En effet, d'après la proposition 2.10, $X(\mathbb{R})$ n'est pas contractible en un point.

### Références bibliographiques